\documentclass[12pt]{amsart}
\usepackage{ifthen,verbatim}

\usepackage{floatflt,graphics}
\usepackage[dvipdfmx]{graphicx}
 
\usepackage{tikz} 
\usepackage{mathrsfs}
\usepackage{amsmath,amssymb,amsthm}

\numberwithin{equation}{section}
\nonstopmode
\theoremstyle{plain}
\newtheorem{cor}[equation]{Corollary}

\newtheorem{lem}[equation]{Lemma}

\newtheorem{thm}[equation]{Theorem}

\theoremstyle{definition}
\newtheorem{rem}[equation]{Remark}

\newtheorem{ex}[equation]{Example}

\newtheorem{defn}[equation]{Definition}

\newenvironment{pf}[1][]{%
 \vskip 3mm
 \noindent
 \ifthenelse{\equal{#1}{}}%
  {{\slshape Proof. }}%
  {{\slshape #1.} }%
 }%
{\qed\bigskip}

\newcounter{alphabet}

\newenvironment{Thm}[1][]{\refstepcounter{alphabet}%
\bigskip%
\noindent%
{\bf Theorem \Alph{alphabet}}%
\ifthenelse{\equal{#1}{}}{}{ (#1)}%
{\bf .}
\itshape}{\vskip 8pt}
\newcounter{minutes}
\divide\time by 60
\newcounter{hours}
\multiply\time by 60
\addtocounter{minutes}{-\time}

\begin{document}
\bibliographystyle{amsplain}
\title[Intrinsic geometry and boundary structure of plane domains]
{
Intrinsic geometry and boundary structure \\ of plane domains 
}

\def\thefootnote{}
\footnotetext{
\texttt{\tiny File:~\jobname .tex,
          printed: \number\year-\number\month-\number\day,
          \thehours.\ifnum\theminutes<10{0}\fi\theminutes}
}
\makeatletter\def\thefootnote{\@arabic\c@footnote}\makeatother

\author[O. Rainio]{Oona Rainio}
\address{Department of Mathematics and Statistics, University of Turku, FI-20014 Turku, Finland}
\email{ormrai@utu.fi}
\author[T. Sugawa]{Toshiyuki Sugawa}
\address{Graduate School of Information Sciences,
Tohoku University, Aoba-ku, Sendai 980-8579, Japan}
\email{sugawa@math.is.tohoku.ac.jp}
\author[M. Vuorinen]{Matti Vuorinen}
\address{Department of Mathematics and Statistics, University of Turku, FI-20014 Turku, Finland}
\email{vuorinen@utu.fi}

\keywords{Condenser capacity, hyperbolic metric, uniformly perfect set.}
\subjclass[2010]{Primary 30F45; Secondary 30C85}
\begin{abstract}
For a non-empty compact set $E$ in a proper subdomain $\Omega$ of the complex plane,
we denote the diameter of $E$ and
the distance from $E$ to the boundary of $\Omega$ by $d(E)$ and $d(E,\partial\Omega),$ respectively.
The quantity $d(E)/d(E,\partial\Omega)$ is invariant under similarities and
plays an important role in Geometric Function Theory.
In the present paper, when $\Omega$ has the hyperbolic distance $h_\Omega(z,w),$ 
we consider the infimum $\kappa(\Omega)$ of
the quantity $h_\Omega(E)/\log(1+d(E)/d(E,\partial\Omega))$ over compact subsets $E$
of $\Omega$ with at least two points, where $h_\Omega(E)$ stands for
the hyperbolic diameter of the set $E.$
We denote the upper half-plane by $\mathbb{H}$.
Our main results claim that $\kappa(\Omega)$ is positive if and only if the boundary
of $\Omega$ is uniformly perfect and that the inequality $\kappa(\Omega)\le \kappa(\mathbb{H})$
holds for all $\Omega,$ where equality holds precisely when $\Omega$ is convex.
\end{abstract}
\thanks{
The authors were supported in part by JSPS KAKENHI Grant Number JP17H02847.
}
\maketitle


\section{Introduction}


Let $\Omega$ be a domain in the complex plane ${\mathbb C}$ with the hyperbolic metric 
$\rho_\Omega(z)|dz|$ of Gaussian curvature $-1$ \cite{bm}. The celebrated 
Uniformization Theorem  \cite[p. 81]{kela} guarantees the existence of $\rho_\Omega$ 
for a domain $\Omega$ when its boundary $\partial \Omega$ contains 
at least three points. Such a domain is called hyperbolic. Here and in what follows, 
the boundary of a domain is taken with respect to the Riemann sphere ${\widehat{\mathbb C}}={\mathbb C}\cup\{\infty\}.$

The function $\rho_\Omega(z)$ is sometimes called the {\it hyperbolic density} of $\Omega.$
For instance, for the unit disk ${\mathbb D}=\{z\in{\mathbb C}: |z|<1\}$ and the upper half-plane
${\mathbb H}=\{z\in{\mathbb C}: {\,\operatorname{Im}\,} z>0\},$ the hyperbolic densities are given by
$\rho_{\mathbb D}(z)=2/(1-|z|^2)$ and $\rho_{\mathbb H}(z)=1/{\,\operatorname{Im}\,} z,$ respectively.
Let $h_\Omega(z_1,z_2)$ denote the {\it hyperbolic distance} induced by 
$\rho_\Omega(z)|dz|$ and $d(z, \partial \Omega)$ the Euclidean distance 
from a point $z\in \Omega$ to the boundary $\partial \Omega.$ Then we have the 
inequality $\rho_\Omega(z)\le 2/d(z, \partial \Omega)$ for each $z\in \Omega$ 
as a simple consequence of Schwarz' Lemma \cite[(2.1)]{BP78}.
On the other hand, the inequality $\rho_\Omega(z)\ge 1/(2d(z, \partial \Omega))$ holds 
for a simply connected domain $\Omega$ 
\cite[(2.2)]{BP78}, \cite[p. 35 Thm 8.6]{bm}, \cite[p.34, (3.2.1)]{gh}. 

The distance on $\Omega$ induced by the continuous Riemannian metric 
$|dz|/d(z, \partial \Omega)$ is called the {\it quasihyperbolic distance} and 
denoted by $k_\Omega(z_1,z_2)$ \cite{GP76}.  We now have the inequality 
$h_\Omega(z_1,z_2)\le 2k_\Omega(z_1,z_2)$ for a general domain $\Omega$
and $h_\Omega(z_1,z_2)\ge k_\Omega(z_1,z_2)/2$ for 
a simply connected domain $\Omega.$ These two inequalities are very handy, 
because there are many estimates for quasihyperbolic distances whereas 
hyperbolic distances are not easy to estimate because the density function 
$\rho_\Omega(z)$ depends on the local boundary structure in the vicinity of $z$ 
in a subtle manner \cite{BP78}, \cite[p.241, Thm 14.5.2]{kela}, \cite{sv}. 
It should be noticed that the second estimate does not apply to general domains, because the hyperbolic distance  is not bounded from below by a constant multiple
of the quasihyperbolic distance, for instance, if the domain has isolated 
boundary points. To measure the similarity between $h_\Omega$ and 
$k_\Omega,$ the domain functional \cite{HM92}
\begin{equation}\label{eq:c}
c(\Omega)=\inf_{z\in \Omega}\rho_\Omega(z)d(z, \partial \Omega)
=\inf_{z_1,z_2\in \Omega,~z_1\ne z_2}\frac{h_\Omega(z_1,z_2)}{k_\Omega(z_1,z_2)}
\end{equation}
is useful, where the second equality will be proven in the next section.
By the above observations, we have $c(\Omega)\le 2$ for a general domain $\Omega$ and
$c(\Omega)\ge 1/2$ for a simply connected domain $\Omega.$
But more is known about this domain constant.

\begin{Thm}
Let $\Omega$ be a hyperbolic domain in ${\mathbb C}.$
Then $c(\Omega)\le 1$ with equality  if and only if $\Omega$ is convex.
Furthermore, $c(\Omega)>0$ if and only if $\partial \Omega$ is uniformly perfect.
\end{Thm}

The general inequality $c(\Omega)\le 1$ is due to Harmelin and Minda \cite{HM92}
and the equality condition is due to Mej\'\i a and Minda \cite{MM90}.
The last assertion is due to Beardon and Pommerenke \cite{BP78}.
Here, a closed set $E$ in ${\widehat{\mathbb C}}$ with ${\operatorname{card}\,}(E)\ge 2$ is said to be 
{\it uniformly perfect} if there is a constant $0<\alpha<1$ such that the closed 
annulus $\alpha r\le |z-a|\le r$ meets $E$ whenever  $a \in E$ and $0<r<{d}(E).$
Here and hereafter, ${\operatorname{card}\,}(E)$ denotes the cardinality of the set $E$ and ${d}(E)$ 
is the Euclidean diameter of $E.$ 
In other words, ${d}(E)=\sup_{z,w\in E} |z-w|.$
We set ${d}(E)=+\infty$ when $\infty\in E.$
For uniformly perfect sets, we refer to \cite{aw}, \cite[pp. 343-345]{gm}, \cite{Pom79}, 
\cite{Pom84}, \cite{SugawaUP} and \cite{SugawaSEUP}.
Uniform perfectness has many applications in potential theory, metric spaces, 
Kleinian groups and complex dynamics  as well as geometric function theory;
see, in addition to the above references, for instance \cite{aw}, \cite{brc}, \cite{SSS} and \cite{WZ}.

In their work about the quasihyperbolic metric, Gehring and Palka \cite{GP76} 
also introduced the {\it distance-ratio metric}
$$
j_\Omega(z_1, z_2)=
\log\left(1+\frac{|z_1-z_2|}{\min\{d(z_1, \partial \Omega),d(z_2, \partial \Omega)\}}\right)
$$
for $z_1,z_2\in \Omega,$ see also \cite[p.61]{HKV20}.
They proved that $j_\Omega(z_1,z_2)\le k_\Omega(z_1,z_2)$ holds always.
It is also known that $j_\Omega$  satisfies the triangle inequality on 
$\Omega$ \cite[p.59, Lemma 4.6]{HKV20}. The opposite inequality characterises 
so called uniform domains: a domain  $\Omega$ is {\it uniform} if and only if 
there exists a constant $b>0$ such that the inequality 
\[k_\Omega(z_1,z_2)\le b j_\Omega(z_1,z_2)\] 
holds,  see  Gehring and Osgood \cite{GO79} and \cite[p.84]{HKV20}. These domains
are ubiquitous in geometric function theory \cite{gh}.

It is a natural and interesting question to ask what can be said if we replace 
$k_\Omega$ by $h_\Omega.$ Our answer is the following result.

\begin{thm}\label{thm:fmt}
Let $\Omega$ be a hyperbolic domain in ${\mathbb C}.$
There is a constant $c>0$ such that $c j_\Omega(z_1,z_2)\le h_\Omega(z_1,z_2)$ for all $z_1,z_2\in\Omega$
if and only if the boundary of $\Omega$ in ${\widehat{\mathbb C}}$ is uniformly perfect.
\end{thm}

In conjunction with the Gehring-Osgood theorem \cite[pp.59-60]{GO79}, 
we have the following result.

\begin{cor}
Let $\Omega$ be a hyperbolic domain in ${\mathbb C}.$ Then the hyperbolic metric 
$h_\Omega$ is comparable with the distance-ratio metric $j_\Omega$
if and only if $\Omega$ is uniform and has uniformly perfect boundary.
\end{cor}

Indeed, if for some constants $0<c_1\le c_2,$ 
\[c_1 j_\Omega(z_1,z_2)\le h_\Omega(z_1,z_2)\le c_2 j_\Omega(z_1,z_2),\quad {\rm for}
\,\,z_1,z_2\in \Omega,\]
 we first see that $\partial \Omega$ is uniformly perfect. Then $h_\Omega$ is 
 comparable with $k_\Omega$ by Theorem A. We now conclude that $\Omega$ 
 is uniform by the Gehring-Osgood theorem. The converse follows readily from 
 Theorem \ref{thm:fmt} and the Gehring-Osgood theorem.

For a subset $E$ of $\Omega$ with ${\operatorname{card}\,}(E)\ge 2,$ we define the set functionals
$$
h_\Omega(E)=\sup_{z_1,z_2\in E}h_\Omega(z_1,z_2)
\quad\text{and}\quad
J_\Omega(E)=\log\left(1+\frac{{d} (E)}{d(E, \partial \Omega)}\right).
$$
Here and hereafter, 
$d(E, F)$ denotes the Euclidean distance between the sets $E$ and $F.$
For a singleton $E=\{z\},$ we write $d(\{z\}, F)=d(z, F)=d(F,z).$
We will use the following monotonicity property frequently in the sequel:
$h_\Omega(E)\le h_\Omega(E')$ and $J_\Omega(E)\le J_\Omega(E')$ for 
$E\subset E'\subset \Omega.$ 
We note that $h_\Omega(E)$ is the hyperbolic diameter of $E$ in $\Omega$ 
and that $J_\Omega(E)$ is important in connection 
with capacity estimates of $E$ (see, for instance, \cite{GSV20}).
We now consider the domain constant
$$
\kappa(\Omega)=\inf_E \frac{h_\Omega(E)}{J_\Omega(E)},
$$
where $E$ runs over all compact subsets of $\Omega$ with ${\operatorname{card}\,}(E)\ge 2.$
As the following result tells, the two domain constants $c(\Omega)$ and 
$\kappa(\Omega)$ are comparable.

\begin{thm}\label{thm:smt}
Let $\Omega$ be a hyperbolic domain in ${\mathbb C}.$ Then the double inequality
$$
\dfrac{c(\Omega)}{2}\le \kappa(\Omega)\le c(\Omega)
$$
holds. In particular, $\kappa(\Omega)>0$ if and only if $\partial \Omega$ 
is uniformly perfect.
\end{thm}

It is a little surprising that the quantity $\kappa(\Omega)$ behaves like 
$c(\Omega)$ in the following sense (compare with Theorem A).

\begin{thm}\label{thm:tmt}
Let $\Omega$ be a hyperbolic domain in ${\mathbb C}.$ Then, the inequality 
$\kappa(\Omega)\le\kappa({\mathbb H})$ holds, where
equality holds if and only if $\Omega$ is convex.
\end{thm}

In view of the above theorem, we are curious about the value of $\kappa({\mathbb H}).$
However, it seems difficult to evaluate it in a simple form.
Since $c({\mathbb H})=1,$ the first part of Theorem \ref{thm:smt} implies
$1/2\le \kappa({\mathbb H})\le 1.$
We will prove later that the inequality $\kappa({\mathbb H})<1$ holds and give
a numerical approximation of the value of $\kappa({\mathbb H})$ in Theorem \ref{thm_approxforOmega}, thus answering 
a problem formulated in \cite[p.455, item (12)]{HKV20}.

The existence of an extremal configuration of the set $E$ for 
the functional $h_\Omega(E)/J_\Omega(E)$ is more subtle.
We will prove the following result in the final section.
We note that a convex domain in ${\mathbb C}$ carries the hyperbolic metric
unless it is ${\mathbb C}$ itself.

\begin{thm}\label{thm:existence}
Let $\Omega$ be a convex proper subdomain of ${\mathbb C}.$ There exists 
a compact subset $E$ in $\Omega$ satisfying $\kappa(\Omega)
=h_\Omega(E)/J_\Omega(E)$ if and only if $\Omega$ is a half-plane.

When $\Omega$ is the upper half-plane ${\mathbb H},$ there exists 
a three-point set $E^*$ of the form $\{i, z_1, z_2\}$ constituting
a hyperbolic equilateral triangle with $\kappa({\mathbb H})=h_{\mathbb H}(E^*)/J_{\mathbb H}(E^*),$
$1<{\,\operatorname{Im}\,} z_j~ (j=1,2)$ and $z_1=-\overline{z_2}.$
Moreover, such an extremal three-point set is unique up to similarities keeping
${\mathbb H}$ invariant.
\end{thm}

In view of the application given in the final section,
it is important to have a lower bound of $\kappa(\Omega)$ when
$\Omega$ is simply connected.
We consider the number
\begin{equation}\label{eq:sc}
\kappa_0=\inf_{\Omega}\kappa(\Omega),
\end{equation}
where $\Omega$ runs over all simply connected proper subdomains of ${\mathbb C}.$
By Theorem \ref{thm:smt} and the well-known estimate $c(\Omega)\ge 1/2,$
we obtain $\kappa_0\ge 1/4.$
On the other hand, when $\Omega$ is the slit domain
$\Omega_0={\mathbb C}\setminus(-\infty,0],$
numerically we have $\kappa(\Omega_0)
\le h_{\Omega_0}(E)/J_{\Omega_0}(E)=0.4251604\dots$
for $E=\{w_0, w_1, w_2\},~ w_0=1,~ 
w_1=2.121820474+1.198476681 i,~ w_2=\bar w_1.$
Note that $h_{\Omega_0}(w_0,w_1)=h_{\Omega_0}(w_0,w_2)
\approx  h_{\Omega_0}(w_1,w_2).$
Thus, we have the following corollary.

\begin{cor}\label{cor}
$1/4\le \kappa_0 < 0.4251605.$
\end{cor}

It is an open problem to determine the value $\kappa_0.$

\medskip

The organization of this paper is as follows.
In Section 2,  preliminary results cocerning the domain constant
$\kappa(\Omega)$ are given
and Theorems \ref{thm:fmt} and \ref{thm:smt} are proved.
Section 3 is devoted to the proof of Theorem \ref{thm:tmt}.
We determine extremal configurations of three-point sets $E$
with respect to the set functional $h_{\mathbb H}(E)/J_{\mathbb H}(E)$
and prove Theorem \ref{thm:existence} in Section 4.
We also give numerical observations on the quantity $\kappa({\mathbb H}).$
We will apply our results to lower estimation of the capacity of a condenser
in the final section.


\section{Preliminaries}


In this section, we prove several simple preliminary results.
We begin with the proof of the second equality in \eqref{eq:c}.
To distinguish the both sides of \eqref{eq:c}, for a while, we write
$$
c(\Omega)=\inf_{z\in \Omega}\rho_\Omega(z)d(z, \partial \Omega)
\quad\text{and}\quad
c'(\Omega)=\inf_{z,w\in \Omega}\frac{h_\Omega(z,w)}{k_\Omega(z,w)}.
$$
We will prove that $c(\Omega)=c'(\Omega).$
Since $\rho_\Omega(z)\ge c(\Omega)/d(z, \partial \Omega),$ 
we easily obtain $h_\Omega(z_1,z_2)\ge c(\Omega) k_\Omega(z_1,z_2).$
Hence, $c'(\Omega)\ge c(\Omega).$
On the other hand, by the formula
$$
\lim_{w\to z}\frac{h_\Omega(z,w)}{k_\Omega(z,w)}
=\lim_{w\to z}\frac{h_\Omega(z,w)}{|z-w|}\cdot\frac{|z-w|}{k_\Omega(z,w)}
=\frac{\rho_\Omega(z)}{1/d(z, \partial \Omega)}=\rho_\Omega(z)d(z, \partial \Omega),
$$
we have $c(\Omega)\ge c'(\Omega).$
Thus, we are done.

\medskip

For the analysis of domain constants, we introduce some variants of the domain constant 
$\kappa(\Omega).$ First, we replace $h_\Omega$ with $k_\Omega$ and define the domain 
constant
$$
\hat\kappa(\Omega)=\inf_E \frac{k_\Omega(E)}{\log(1+{d}(E)/d(E, \partial \Omega))}\,,
$$
where the infimum is taken over all compact subsets $E$ of $\Omega$ with ${\operatorname{card}\,}(E)\ge2.$
Here, $k_\Omega(E)$ denotes the quasihyperbolic diameter of $E.$
We also define the following auxiliary domain constants for integers $n\ge 2$:
$$
\kappa_n(\Omega)=\inf_{E\subset \Omega, {\operatorname{card}\,}(E)=n} 
\frac{h_\Omega(E)}{\log(1+{d}(E)/d(E, \partial \Omega))}
$$
and
$$
\hat\kappa_n(\Omega)=\inf_{E\subset \Omega, {\operatorname{card}\,}(E)=n} 
\frac{k_\Omega(E)}{\log(1+{d}(E)/d(E, \partial \Omega))}\,.
$$
For $E=\{z_1,\dots, z_n\},$ letting $z_n\to z_{n-1},$ we observe that
$$
\kappa_2(\Omega)\ge \kappa_3(\Omega)\ge \cdots \ge \kappa(\Omega)
$$
and
$$
\hat\kappa_2(\Omega)\ge \hat\kappa_3(\Omega)\ge \cdots \ge \hat\kappa(\Omega)\,.
$$
For these domain constants, we have the following results.
In particular, we see that $\kappa_n(\Omega)=\kappa(\Omega)$ and 
$\hat\kappa_n(\Omega)=\hat\kappa(\Omega)$ for every $n\ge 3.$

\begin{lem}\label{lem:gamma_n}
\begin{enumerate}
\item[(i)]
$\hat\kappa_2(\Omega)\ge 1\,.$
\vspace{1mm}
\item[(ii)]
$\kappa_3(\Omega)=\kappa(\Omega)$ and $\hat\kappa_3(\Omega)=\hat\kappa(\Omega)\,.$
\vspace{1mm}
\item[(iii)]
$\kappa_2(\Omega)\le 2\kappa_3(\Omega)$ and $\hat\kappa_2(\Omega)\le 2\hat\kappa_3(\Omega)\,.$
\end{enumerate}
\end{lem}

\begin{pf}
Part (i) is clear from the Gehring-Palka inequality $k_\Omega(z_1, z_2)\ge j_\Omega(z_1, z_2).$
Let $E$ be an arbitrary compact set in $\Omega$ with ${\operatorname{card}\,}(E)\ge 2.$
Take $z_0, z_1,z_2\in E$ so that ${d}(E)=|z_1-z_2|$ 
and $d(E, \partial \Omega)=d(z_0, \partial \Omega)$
and let $E_0=\{z_0, z_1,z_2\}.$
(Note that one of the points $z_1, z_2$ may be the same as $z_0.$) Then
\begin{align}
\notag
h_\Omega(E)\ge h_\Omega(E_0) 
&\ge\kappa_3(\Omega)\log(1+{d}(E_0)/d(E_0, \partial \Omega)) \\
\label{eq:chain}
&= \kappa_3(\Omega)\log(1+|z_1-z_2|/d(z_0, \partial \Omega)) \\
\notag
&=\kappa_3(\Omega)\log(1+{d}(E)/d(E, \partial \Omega)).
\end{align}
Taking the infimum over compact subsets $E$ of $\Omega,$
we obtain the inequality $\kappa(\Omega)\ge\kappa_3(\Omega).$
Since $\kappa(\Omega)\le\kappa_3(\Omega)$ as we noted above, 
we conclude $\kappa(\Omega)=\kappa_3(\Omega).$
In the same way, we can verify $\hat\kappa(\Omega)=\hat\kappa_3(\Omega).$

Finally, we prove part (iii).
Let $E\subset \Omega$ with ${\operatorname{card}\,}(E)=3$ and choose $z_0\in E$ 
so that $d(E, \partial \Omega)=d(z_0, \partial \Omega).$
Also choose $z_1,z_2\in E$ so that ${d}(E)=|z_1-z_2|.$
Then
\begin{align*}
& \, \quad  \log(1+{d}(E)/d(E, \partial \Omega)) \\
&=\log(1+|z_1-z_2|/d(z_0, \partial \Omega)) \\
&\le\log(1+(|z_1-z_0|+|z_2-z_0|)/d(z_0, \partial \Omega)) \\
&\le\log(1+|z_1-z_0|/d(z_0, \partial \Omega))+\log(1+|z_2-z_0|/d(z_0, \partial \Omega)) \\
&\le \kappa_2(\Omega)^{-1}(h_\Omega(z_1,z_0)+h_\Omega(z_2,z_0)) \\
&\le 2h_\Omega(E)/\kappa_2(\Omega),
\end{align*}
which implies $\kappa_2(\Omega)\le 2\kappa_3(\Omega).$
In the same way, we can prove the other inequality.
\end{pf}

We need also the following simple lemma.

\begin{lem}\label{lem:g2}
For a hyperbolic domain $\Omega$ in ${\mathbb C},$ the inequality
$\kappa_2(\Omega)\le c(\Omega)$ holds.
\end{lem}

\begin{pf}
Noting the formula
$$
\lim_{w\to z}\frac{h_\Omega(z,w)}{j_\Omega(z,w)}=\rho_\Omega(z)d(z,\partial \Omega),
$$
we have
$$
\kappa_2(\Omega)=\inf_{z\ne w}\frac{h_\Omega(z,w)}{j_\Omega(z,w)}
\le\inf_{z\ne w}\rho_\Omega(z)d(z,\partial \Omega)=c(\Omega).
$$
\end{pf}

We are now in a position to prove Theorem \ref{thm:smt}.

\begin{pf}[Proof of Theorem \ref{thm:smt}]
By the above lemma and the inequality $h_\Omega(x,y)\ge c(\Omega) k_\Omega(x,y),$
for an arbitrary compact set $E$ in $\Omega,$
we have
\begin{align*}
\frac{h_\Omega(E)}{\log(1+{d}(E)/d(E, \partial \Omega))}
&\ge \frac{c(\Omega)k_\Omega(E)}{\log(1+{d}(E)/d(E, \partial \Omega))} \\[2mm]
&\ge c(\Omega)\hat\kappa(\Omega)
\ge \frac{c(\Omega)}2\hat\kappa_2(\Omega) \ge \frac{c(\Omega)}2.
\end{align*}
Hence we have $\kappa(\Omega)\ge c(\Omega)/2.$
The other inequality follows from Lemma \ref{lem:g2}:
$$
\kappa(\Omega)\le \kappa_2(\Omega)\le c(\Omega).
$$
\end{pf}

We now prove Theorem \ref{thm:fmt}.

\begin{pf}[Proof of Theorem \ref{thm:fmt}]
Assume that $c\,j_\Omega(z_1,z_2)\le h_\Omega(z_1,z_2)$ for $z_1,z_2\in \Omega.$
Then $\kappa_2(\Omega)\ge c.$
By Lemma \ref{lem:gamma_n} and Theorem \ref{thm:smt}, we obtain
$$
c(\Omega)\ge \kappa(\Omega)\ge \frac12\kappa_2(\Omega)\ge \frac{c}2>0.
$$
Thus, $\partial \Omega$ is uniformly perfect. Conversely, if $\partial \Omega$ 
is uniformly perfect, similarly we obtain
$\kappa_2(\Omega)\ge\kappa(\Omega)\ge c(\Omega)/2>0.$
Thus, $c\,j_\Omega(z_1,z_2)\le h_\Omega(z_1,z_2)$ holds with $c=\kappa_2(\Omega)>0.$
\end{pf}


\section{Proof of Theorem \ref{thm:tmt}}\label{sec:ext}


In this section, we will prove Theorem \ref{thm:tmt} step by step.
We begin with the following result.

\begin{lem}\label{lem:D}
For any hyperbolic domain $\Omega$ in ${\mathbb C},$ 
the inequality $\kappa(\Omega)\le \kappa({\mathbb D})$ holds.
\end{lem}

\begin{pf}
By definition, for a given $\varepsilon>0,$ 
there is a compact subset $E$ of ${\mathbb D}$ such that
$$
\frac{h_{\mathbb D}(E)}{J_{\mathbb D}(E)}<\kappa({\mathbb D})+\varepsilon.
$$
Moreover, by rotating $E$ if necessary, we may further assume that
the nearest point of the boundary $\partial{\mathbb D}$ to $E$ is $1.$
Namely, $d(E, \partial{\mathbb D})={d}(E, 1).$

Let $\Omega$ be an arbitrary hyperbolic domain in ${\mathbb C}.$ For an arbitrarily 
fixed point $z_0\in\Omega,$ choose $\zeta_0\in\partial\Omega$
so that $d(z_0,\partial\Omega)=|z_0-\zeta_0|.$
Since $\kappa(\Omega)$ is invariant under similarities,
we may assume that $z_0=0$ and $\zeta_0=1.$
Then ${\mathbb D}\subset \Omega.$
By the domain monotonicity of the hyperbolic metric, we have 
$h_\Omega(E)\le h_{\mathbb D}(E).$
On the other hand, we have $d(E, \partial\Omega)={d}(E,1)=d(E,\partial{\mathbb D})$
so that $J_\Omega(E)=J_{\mathbb D}(E).$ Hence,
$$
\kappa({\mathbb D})+\varepsilon>\frac{h_{\mathbb D}(E)}{J_{\mathbb D}(E)}
\ge\frac{h_\Omega(E)}{J_\Omega(E)}\ge \kappa(\Omega).
$$
Since $\varepsilon>0$ is arbitrary, we obtain the required inequality 
$\kappa({\mathbb D})\ge\kappa(\Omega).$
\end{pf}

\begin{rem}\label{onjdiam} {\rm Note that the set functional $J_D(E)$ in the above proof is not the same thing as the diameter of $E$ in the $j_D$ metric
\[
j_D(E) = \sup\{j_D(x,y):  x,y \in E\}\,.
\]
It is easy to see that the inequality
\[
J_D(E)/2 \le j_D(E) \le J_D(E)
\]
holds for all  $E \subset D$, with equality in the second 
inequality if $E$ is a disk, ${\operatorname{card}\,}(E)=2$, or ${\operatorname{card}\,}(E)=3$ and the triangle with vertices $E$ is either equilateral or a so-called
Reuleaux triangle.
 }

\end{rem}

Moreover, for a half-plane, we have the following result.

\begin{lem}\label{lem:H}
Let $H$ be an open half-plane in ${\mathbb C}.$ Then $\kappa({\mathbb D})=\kappa(H).$
\end{lem}

\begin{pf}
By Lemma \ref{lem:D}, it is enough to prove the inequality $\kappa(H)\ge\kappa({\mathbb D}).$
We choose the right half-plane $\{z\,:\,{\,\operatorname{Re}\,} z>0\}$ as $H.$
For every $\varepsilon>0,$ we can find a compact subset $E$ of $H$ such that
$$
\frac{h_H(E)}{J_H(E)}<\kappa(H)+\varepsilon.
$$
Let $\zeta_0$ be the nearest boundary point to $E.$ For simplicity, we assume 
that $\zeta_0=0.$ For $R>0,$ we denote the disk 
$\{z \,:\, |z-R|<R\}$ by $\Delta_R$. For a large enough  $R,$ $E\subset \Delta_R$ 
and $d(E, \partial\Delta_R)={d}(E,0)=d(E,\partial H)$ so that
$J_H(E)=J_{\Delta_R}(E).$ On the other hand, since 
$$
\rho_{\Delta_R}(z)=\frac{2R}{R^2-|z-R|^2}
=\frac1{{\,\operatorname{Re}\,} z-|z|^2/(2R)} \to \frac1{{\,\operatorname{Re}\,} z}=\rho_H(z)
$$
locally uniformly on ${\mathbb H},$ we obtain $h_{\Delta_R}(E)\to h_H(E)$ as $R\to+\infty.$
Noting the inequality 
$$h_{\Delta_R}(E)/J_{\Delta_R}(E)\ge \kappa(\Delta_R)=\kappa({\mathbb D}),$$
we have
$$
\frac{h_H(E)}{J_H(E)}
=\lim_{R\to+\infty}\frac{h_{\Delta_R}(E)}{J_{\Delta_R}(E)}\ge \kappa({\mathbb D}).
$$
Hence, $\kappa(H)+\varepsilon>\kappa({\mathbb D}).$ Since $\varepsilon>0$ was arbitrary, 
we obtain the inequality $\kappa(H)\ge\kappa({\mathbb D})$ as required.
\end{pf}

We next prove the following lemma.

\begin{lem}\label{lem:conv}
Let $\Omega$ be a convex domain in ${\mathbb C}$ with $\Omega\ne{\mathbb C}.$
Then $\kappa(\Omega)=\kappa({\mathbb D}).$
\end{lem}

\begin{pf}
Let $E$ be any compact subset of $\Omega.$
Take $\zeta_0\in\partial\Omega$ so that $d(E,\partial \Omega)={d}(E,\zeta_0).$
Since $\Omega$ is convex, there is a supporting line, say, $L$ at the point $\zeta_0.$
Let $H$ be the connected component of ${\mathbb C}\setminus L$ containing $\Omega.$
Then $\Omega\subset H$ and $\zeta_0\in\partial H=L.$
Since $d(E,\partial H)={d}(E,\zeta_0)=d(E,\partial \Omega),$ we obtain
$$
\frac{h_\Omega(E)}{J_\Omega(E)}
\ge \frac{h_H(E)}{J_H(E)}\ge\kappa(H)=\kappa({\mathbb D}).
$$
Here, we used Lemma \ref{lem:H}.
Taking the infimum over $E,$ we obtain the inequality $\kappa(\Omega)\ge \kappa({\mathbb D}).$
Recalling Lemma \ref{lem:D}, we have the desired relation.
\end{pf}

To deduce the equality condition is the most subtle part in the proof of
Theorem \ref{thm:tmt}.
A key ingredient is Keogh's lemma about non-convex domains.
See Figure 1.



\begin{figure}[ht!]
    \centering
    \begin{tikzpicture}[scale=3]
    \draw[line width=0.35mm] (0,0) circle (1cm);
    \draw[line width=0.35mm] (0.839,1.149) arc (130:230:1.5);
    \draw[line width=0.35mm] (0.303,0) circle (0.03cm);
    \draw[line width=0.35mm, xshift=0.3cm] (0.3,-1.3) .. controls (0.5,-0.7)  and (0.5,-0.6) .. (0,0) .. controls (0.3,0.5) and (0.4,0.6) .. (0.1,1.5);
    \node[scale=1.3] at (0,0.4) {$G$};
    \node[scale=1.3] at (0.25,1.3) {$\Omega$};
    \node[scale=1.3] at (-0.4,0.75) {$\Delta_1$};
    \node[scale=1.3] at (0.9,1) {$\Delta_2$};
    \node[scale=1.3] at (0.45,0.03) {$\zeta_0$};
    \end{tikzpicture}
    \label{fig0}
    \caption{The domain $G=\Delta_1\backslash\overline{\Delta_2}$ in $\Omega$}
\end{figure}

\begin{lem}[Keogh \cite{Keogh76}]
Suppose that a domain $\Omega$ in ${\mathbb C}$ is not convex.
Then there are two open disks $\Delta_1$ and $\Delta_2$ whose boundaries
intersect perpendicularly such that $G=\Delta_1\setminus \overline{\Delta_2}$
is contained in $\Omega$ and the midpoint $\zeta_0$ of the concave boundary arc
$\Delta_1\cap\partial\Delta_2$ of $G$ lies on the boundary $\partial\Omega$ of $\Omega.$
\end{lem}

We are now ready to prove the following result, which is the last piece of
the proof of Theorem \ref{thm:tmt}.

\begin{lem}\label{lem:non-conv}
Let $\Omega$ be a non-convex domain in ${\mathbb C}.$
Then $\kappa(\Omega)<\kappa({\mathbb D}).$
\end{lem}

\begin{pf}
We find open disks $\Delta_1, \Delta_2$ as in Keogh's lemma
so that $G=\Delta_1\setminus \overline{\Delta_2}\subset\Omega$
and the midpoint $\zeta_0$ of the concave boundary arc of $G$ 
is contained in  $\partial\Omega.$
We may assume that $\Delta_1={\mathbb D}$ and $\zeta_0=a\in(0,1)$ so that
the center of $\Delta_2$ lies on the real axis.
Then the second disk $\Delta_2$ is the image of the right half-plane $H$ under the
M\"obius transformation
$$
T(z)=\frac{z+a}{1+az}.
$$
Thus, $G=T({\mathbb D}_-),$ where ${\mathbb D}_-$ is the left half $\{z\in{\mathbb D}: {\,\operatorname{Re}\,} z<0\}$ of the unit disk.
We now construct a conformal map $f$ of the upper half-plane ${\mathbb H}$ onto $G$ as follows.
We denote the analytic automorphism $(1+z/2)/(1-z/2)$ of ${\mathbb H}$ by $M$.
Note that $M$ maps the positive imaginary axis $i\,{\mathbb R}_+=\{iy: 0<y<+\infty\}$
onto the upper half of the unit circle $|\zeta|=1.$
The function $S(\zeta)=\sqrt{\zeta}$ maps ${\mathbb H}$ onto the first quadrant 
$D=\{w: {\,\operatorname{Re}\,} w>0, {\,\operatorname{Im}\,} w>0\}.$
Then the M\"obius transformation $L(w)=i(w-1)/(w+1)$ maps $D$ onto
the left half ${\mathbb D}_-$ of ${\mathbb D}.$
Hence, the function $f=T\circ L\circ S\circ M$ maps ${\mathbb H}$ onto $G$
in such a way that $f(i\,{\mathbb R}_+)=(-1,a).$
More concretely, $f$ is expressed by
$$
f(z)=T\left(i\frac{\sqrt{1+z/2}-\sqrt{1-z/2}}{\sqrt{1+z/2}+\sqrt{1-z/2}}\right).
$$
In view of this form, we see that $f(z)$ is analytic on $|z|<1.$
(This follows also from the Schwarz reflection principle.)
Therefore, we can expand $f(z)$ about $z=0$ as follows:
$$
f(z)=a+a_1z+a_2z^2+\cdots \quad (|z|<1).
$$
By a straightforward computation, we have here
$$
a_1=\frac i4(1-a^2), \quad
a_2=\frac{a}{16}(1-a^2)
$$
and therefore
\begin{equation}\label{eq:a12}
A:=\frac{a_2}{a_1}=\frac{a}{4i}.
\end{equation}

Let $E_x:=xE^*=\{xz_j: j=0,1,2\}$ for $0<x<1,$
where $E^*=\{z_0, z_1, z_2\}\subset{\mathbb H}$ with $z_0=i$ is the set in 
Theorem \ref{thm:existence} and thus $\kappa({\mathbb H})=h_{\mathbb H}(E^*)/J_{\mathbb H}(E^*).$
Let $w_j=f(xz_j)$ and set $E_x'=f(E_x)=\{w_j: j=0,1,2\}.$
Since $f(xz)=a+a_1xz+O(x^2)$ as $x\to0$ locally uniformly in $z,$
${d}(E_x')=|w_1-w_2|$ and 
$d(E_x', \partial G)=d(w_0,\partial G)=d(w_0,\Delta_1\cap\partial\Delta_2)$
for a small enough $x>0.$
Note here that $w_0=f(xz_0)=f(ix)\in (0,a)$ because $f(i\,{\mathbb R}_+)=(-1,a).$
Hence, $d(E_x', \partial G)=d(w_0,\Delta_1\cap\partial\Delta_2)=d(w_0,a).$
We now look at the quantity
$$
F(x)=\frac{{d}(E_x')}{d(E_x', \partial G)}
=\frac{|w_1-w_2|}{|w_0-a|}
=\left|\frac{w_1-w_2}{w_0-a}\right|.
$$
We observe that
$$
W=\frac{w_1-w_2}{w_0-a}
=\frac{f(xz_1)-f(xz_2)}{f(xz_0)-f(0)}
$$
is even analytic in $x\in{\mathbb D}$ and we compute
\begin{align*}
W&=\frac{a_1x(z_1-z_2)+a_2x^2(z_1^2-z_2^2)+O(x^3)}{a_1xz_0+a_2x^2z_0^2+O(x^3)} \\
&=\frac{z_1-z_2}{z_0}\cdot \frac{1+Ax(z_1+z_2)+O(x^2)}{1+Axz_0+O(x^2)} \\
&=\frac{z_1-z_2}{z_0}\cdot \big[1+Ax(z_1+z_2-z_0)+O(x^2)\big],
\end{align*}
where $A=a_2/a_1=a/(4i)$ by \eqref{eq:a12}.
Hence $F(x)=|W|$ is real analytic in $-1<x<1$ and
\begin{align*}
F(x)
&=\frac{|z_1-z_2|}{|z_0|}\Big\{1+{\,\operatorname{Re}\,}\big[ Ax(z_1+z_2-z_0)\big]+O(x^2)\Big\} \\
&=\frac{|z_1-z_2|}{|z_0|}\left\{1+\frac{ax}{4}{\,\operatorname{Im}\,}(z_1+z_2-z_0)+O(x^2)\right\}
\end{align*}
as $x\to0.$
Since ${\,\operatorname{Im}\,} z_j={d}(z_j,\partial{\mathbb H})>{d}(z_0,\partial{\mathbb H})$ for $j=1,2,$
we have
$$
F(0)=\frac{|z_1-z_2|}{|z_0|}=\frac{{d}(E^*)}{d(E^*,\partial{\mathbb H})}
\quad\text{and}\quad
F'(0)=\frac{a|z_1-z_2|}{4|z_0|}{\,\operatorname{Im}\,}(z_1+z_2-z_0)>0.
$$
In particular, $F(x)$ is strictly increasing at $x=0$ and thus
$F(x)>F(0)$ for small enough $x>0.$
Since $G\subset\Omega,$ we have the inequality $h_\Omega(E_x')\le h_G(E_x').$
We also note that 
$$
d(E_x',\partial\Omega)\ge d(E_x', \partial G)=d(w_0,a)\ge d(E_x',\partial\Omega),
$$
because $a\in\partial\Omega,$ and therefore $d(E_x',\partial\Omega)=d(E_x', \partial G)$
so that $J_\Omega(E_x')=J_G(E_x').$
Moreover, since the hyperbolic distance is conformally invariant, $h_G(E_x')
=h_G(f(E_x))=h_{\mathbb H}(E_x)=h_{\mathbb H}(E^*).$
Hence, for a small enough $x>0,$
\begin{align*}
\kappa(\Omega)
&\le\frac{h_\Omega(E_x')}{J_\Omega(E_x')} 
\le\frac{h_G(E_x')}{J_G(E_x')} 
=\frac{h_{\mathbb H}(E^*)}{\log(1+F(x))} \\
&<\frac{h_{\mathbb H}(E^*)}{\log(1+F(0))} 
=\frac{h_{\mathbb H}(E^*)}{J_{\mathbb H}(E^*)}=\kappa({\mathbb H}).
\end{align*}
The proof is finished.
\end{pf}

Now Theorem \ref{thm:tmt} follows from Lemmas \ref{lem:D},
\ref{lem:conv} and \ref{lem:non-conv}.


\section{Extremal configuration of three points in ${\mathbb H}$}


In this section, we work to find extremal configurations of three-point sets $E$ in the
upper half-plane for the functional $h_{\mathbb H}(E)/J_{\mathbb H}(E).$
Since the both quantities $h_{\mathbb H}(E)$ and $J_{\mathbb H}(E)$ are invariant under
the affine mappings of the form $z\mapsto az+b$ with $a>0$ and $b\in{\mathbb R},$
we may restrict our attention to the family ${\mathcal E}$ of three-point subsets $E$ of ${\mathbb H}$
containing $i=\sqrt{-1}$ with $d(E,\partial{\mathbb H})=d(i,\partial{\mathbb H})=1.$
Namely, the infimum in the definition of $\kappa_3({\mathbb H})$ may be limited to ${\mathcal E}:$
$$
\kappa_3({\mathbb H})
=\inf_{E\in{\mathcal E}}\frac{h_{\mathbb H}(E)}{J_{\mathbb H}(E)}
=\inf_{E\in{\mathcal E}}\frac{h_{\mathbb H}(E)}{\log(1+{d}(E))}.
$$
Our goal in this section is to determine the extremal sets $E$ for which the above
infimum is attained, and to compute (at least numerically) the value of $\kappa_3({\mathbb H}).$
First, we note the following fact for the upper half-plane ${\mathbb H}.$
Though the result is essentially known (e.g., \cite[Lemma 4.9 (2)]{HKV20}),
we give a short proof for convenience of the reader.

\begin{lem}\label{lem:gm1}
$$
\kappa_2({\mathbb H})=\inf_{z_1,z_2\in{\mathbb H}}
\frac{h_{\mathbb H}(z_1,z_2)}{j_{\mathbb H}(z_1,z_2)}=1.
$$
\end{lem}

\begin{pf}
Note that $\rho_{\mathbb H}(z)=1/{\,\operatorname{Re}\,} z=1/{d}(z,\partial{\mathbb H}).$
Hence, we have $h_{\mathbb H}(z,w)=k_{\mathbb H}(z,w)$ for $z,w\in{\mathbb H}.$
Thus, the inequality $j_{\mathbb H}(z,w)\le h_{\mathbb H}(z,w)$ is nothing but the Gehring-Palka inequality \cite{GP76}.
Hence, we have $\kappa_2({\mathbb H})\ge 1.$
On the other hand, by Lemma \ref{lem:g2}, we have
$\kappa_2({\mathbb H})\le c({\mathbb H})\le 1,$ where the last inequality follows from Theorem A.
\end{pf}

We will write
$$
\Delta(z_0, r)=\{z\in{\mathbb H}: h_{\mathbb H}(z,z_0)< r\}
=\{z: |z-z_0|< \rho|z-\bar z_0|\}
$$
for the open hyperbolic disk in ${\mathbb H}$ centered at $z_0\in{\mathbb H}$ 
with hyperbolic radius $r>0,$
where $\rho=\tanh(r/2)=(e^r-1)/(e^r+1)\in(0,1)$ 
and denote its closure by ${\overline{\Delta}}(z_0, r)$.
We need the following elementary fact for the proof of Lemma \ref{lem:extremal},
which will be a key result below.

\begin{lem}\label{lem:circle}
Let $C$ be the boundary circle of the hyperbolic disk $\Delta(z_0,r)$ in ${\mathbb H}.$
\begin{enumerate}\item[(i)]
The Euclidean distance $|z-z_0|$ between $z\in C$ and $z_0$
takes its maximum at the top of $C$ and its minimum at the bottom of $C.$
\vspace{1mm}
\item[(ii)]
The Euclidean diameter of the circle $C$ is $2({\,\operatorname{Im}\,} z_0)\sinh r.$
\vspace{1mm}
\item[(iii)]
The hyperbolic distance of the endpoints of an arbitrary diameter of the circle $C$
is at least equal to $\varphi(r)$ given in \eqref{eq:phi}.
\end{enumerate}
\end{lem}

\begin{pf}
We write $z_0=x_0+iy_0.$
It is well known (see, e.g., \cite[(4.11)]{HKV20}) 
that the boundary of $\Delta(z_0,r)$ is the Euclidean circle
$|z-c|=R,$ where
$$
c=x_0+iy_0 \cosh r
\quad\text{and}\quad
R=y_0 \sinh r.
$$
Since ${\,\operatorname{Re}\,} z_0={\,\operatorname{Re}\,} c$ and ${\,\operatorname{Im}\,} z_0<{\,\operatorname{Im}\,} c$, it is evident that
$|z-z_0|$ is maximized at $z=c+iR$ and minimized at $z=c-iR$ on $C.$
The proof of the first assertion is now complete.
The second assertion is clear because the Euclidean diameter of $C$ is $2R.$
It is clear that the diameter of the circle $C$ with the minimal
hyperbolic diameter is $[c-R,c+R].$
We now compute the hyperbolic distance
\begin{align}
\notag
h_{\mathbb H}(c+R,c-R)&=h_{\mathbb H}(i\cosh r+\sinh r, i\cosh r-\sinh r) \\
\notag
&=2\,\mathrm{artanh}\, \frac{\sinh r}{\sqrt{\cosh 2r}} \\
\notag
&=\log\frac{\sqrt{\cosh 2r}+\sinh r}{\sqrt{\cosh 2r}-\sinh r} \\
&=2\log\frac{\sqrt{\cosh 2r}+\sinh r}{\cosh r}
=:\varphi(r).
\label{eq:phi}
\end{align}
Then the third assertion follows.
\end{pf}

\begin{rem}\label{rem:circle}
By geometry, we see that $|c+i\, Re^{\pm i\theta}-z_0|$ is strictly
decreasing in $0<\theta<\pi,$ which will be needed in the proof of Lemma
\ref{lem:extremal}.

We remark also that the sharp upper bound of the hyperbolic distance
of the endpoints of a diameter of $C$ is $h_{\mathbb H}(c+iR, c-iR)=2r.$
By the form of $\varphi(r),$ we also see that
$\varphi(r)\to \log\frac{\sqrt 2+1}{\sqrt 2-1}=2\log(\sqrt 2+1)
=1.7627\dots$ as $r\to+\infty.$
\end{rem}

In order to find the extremal configuration, we divide the family ${\mathcal E}$
into one-parameter subfamilies.
More concretely, for $u>0,$ let ${\mathcal E}(u)$ be the subfamily of ${\mathcal E}$
consisting of sets $E$ with $h_{\mathbb H}(E)=2u.$
Then
\begin{equation}\label{eq:g3}
\kappa_3({\mathbb H})
=\inf_{0<u<+\infty}\inf_{E\in{\mathcal E}(u)} \frac{2u}{J_{\mathbb H}(E)}
=\inf_{0<u<+\infty}\frac{2u}{\log(1+M(u))},
\end{equation}
where
\begin{equation}\label{eq:M}
M(u)=\sup_{E\in{\mathcal E}(u)}{d}(E)
\end{equation}
Our task is to find the extremal configuration
of $E\in{\mathcal E}(u)$ for the functional ${d}(E).$
We first define a candidate of the extremal set.
For a given number $u>0,$ we choose $t>0$ and $\theta\in(0,\pi/2)$
such that
$$
h_{\mathbb H}(ie^{t+i\theta},ie^{t-i\theta})=h_{\mathbb H}(ie^{t+i\theta},i)=2u.
$$
In other words, we choose $t$ and $\theta$ so that the set
$E^*(u)=\{i, ie^{t+i\theta}, ie^{t-i\theta}\}$
forms the vertices of a hyperbolic equilateral triangle with sidelength $2u.$
We now give formulae describing $\theta$ and $t$ in terms of $u.$
Since $h_{\mathbb H}(ie^{t+i\theta},ie^t)=u,$ we obtain 
$u=2\,\mathrm{artanh}\,(\tan(\theta/2))$ and thus 
\begin{equation}\label{eq:theta}
\theta=2\arctan(\tanh(u/2)).
\end{equation}
Moreover, by the hyperbolic cosine formula for a hyperbolic right triangle 
\cite[Thm 7.11.1, p. 146]{be}, we have
$$
\cosh t=\cosh h_{\mathbb H}(ie^t,i)
=\frac{\cosh h_{\mathbb H}(ie^{t+i\theta},i)}{\cosh h_{\mathbb H}(ie^{t+i\theta},ie^t)}
=\frac{\cosh 2u}{\cosh u}.
$$
Hence, 
\begin{equation}\label{eq:t}
t=\,\mathrm{arcosh}\,((\cosh 2u)/\cosh u).
\end{equation}
We now compute
$$
|ie^{t+i\theta}-ie^{t-i\theta}|=2e^t\sin\theta=\chi(u),
$$
where
\begin{align}\label{eq:chi}
\chi(u)&=2e^{\,\mathrm{arcosh}\,((\cosh 2u)/\cosh u)}\sin \big[2\arctan\tanh(u/2)\big] \\
\notag
&=2\frac{\cosh 2u+\sqrt{(\cosh^22u)-(\cosh^2u)}}{\cosh u}\cdot\tanh u \\
\notag
&=\frac{2\sinh u}{1+\sinh^2u}\big[1+2\sinh^2u+\sinh u\sqrt{3+4\sinh^2u}\big].
\end{align}
Note that $\chi(u)\le{d}(E^*(u)).$
In the same way, we compute
\begin{align*}
{\,\operatorname{Im}\,} \big(ie^{t+i\theta}\big)&=e^t\cos\theta
=e^{\,\mathrm{arcosh}\,((\cosh 2u)/\cosh u)}\cos \big[2\arctan\tanh(u/2)\big] \\
&=\frac{\cosh 2u+\sqrt{\cosh^22u-\cosh^2u}}{\cosh u}\cdot\frac1{\cosh u} 
>\frac{\cosh 2u}{\cosh^2u}>1\,.
\end{align*}
Therefore, we obtain $d(E^*(u),\partial{\mathbb H})=1$ for every $u>0.$
We summarize the above observations in the following lemma.

\begin{lem}\label{lem:E}
The set $E^*(u)$ of the vertices of the hyperbolic equilateral triangle in ${\mathbb H}$
with sidelength $2u$ constructed above belongs to ${\mathcal E}(u)$ for every $u>0.$
\end{lem}

We make further preparatory observations.

\begin{lem}\label{lem:11/4}
If $0<u\le \log(11/4)\approx 1.0116$, then ${d}(E^*(u))=\chi(u)$ and
$$
\frac{2u}{\log(1+M(u))}<1.
$$
\end{lem}

\begin{pf}
We will prove the inequality
\begin{equation}\label{eq:<1}
\frac{2u}{\log(1+\chi(u))}<1
\end{equation}
for $0<u\le \log(11/4).$
Since $E^*(u)\in{\mathcal E}(u)$ by Lemma \ref{lem:E}, we have
$M(u)\ge {d} (E^*(u))\ge \chi(u).$
Thus, the second assertion will follow from \eqref{eq:<1}.

By using the elementary inequality $\sqrt{3+4\sinh^2u}>\sqrt{3+3\sinh^2u}
=\sqrt3\cosh u$ for $u>0,$
we obtain the estimate
$$
\chi(u)
>\frac{2\sinh u}{1+\sinh^2u}\big[1+2\sinh^2u+\sqrt 3\sinh u\cosh u\big].
$$
Thus, we have
\begin{align*}
\chi(u)+1-e^{2u}
&\ge \frac{2\sinh u}{1+\sinh^2u}\big[1+2\sinh^2u+\sqrt 3\sinh u\cosh u\big]+1-e^{2u} \\
&=\frac{(\sqrt3-1)(e^u+1)(e^u-1)^2P(e^u-1)}{e^u(e^{2u}+1)^2},
\end{align*}
where $P(T)$ is the polynomial given by
$$
P(T)=4+(7+\sqrt 3)T+4T^2-\sqrt 3T^3-\frac{1+\sqrt 3}2 T^4.
$$
We now estimate $P(T)$ for $T\ge0$ from below:
$$
P(T)\ge 4+8T+4T^2-2T^3-2T^4=2(1+T)(2+2T-T^3).
$$
Since $Q(T)=2+2T-T^3$ is concave on $[0,+\infty),$ we have
$$Q(T)\ge \min\{Q(0),Q(7/4)\}=9/64>0\quad {\rm for}\,\, 0\le T\le 7/4.$$
Hence, we have proved that $e^{2u}<1+\chi(u)$ and
thus \eqref{eq:<1} holds for $0<u\le \log(11/4).$

Finally, we prove that ${d}(E^*(u))=\chi(u)$ for such $u.$
Indeed, the inequality 
\[|ie^{t+i\theta}-ie^{t-i\theta}|<|ie^{t+i\theta}-i|\]
would hold otherwise. Then the two-point subset $E=\{i, ie^{t+i\theta}\}$ 
of $E^*(u)$ satisfies $h_{\mathbb H}(E)=2u,~{d}(E)={d}(E^*(u))$ and 
$d(E,\partial{\mathbb H})=d(E^*(u),\partial{\mathbb H})=1.$
Thus, we would have
$$
\frac{2u}{\log(1+\chi(u))}>
\frac{2u}{J_{\mathbb H}(E^*(u))}=\frac{2u}{J_{\mathbb H}(E)}\ge \kappa_2({\mathbb H})=1
$$
by Lemma \ref{lem:gm1}.
This contradicts \eqref{eq:<1}.
In this way, we have proved that ${d}(E^*(u))=\chi(u).$
\end{pf}

\begin{lem}\label{lem:u0}
Let $0<u<+\infty.$
The condition $\varphi(2u)\ge 2u$ holds if and only if $u\le u_0,$
where $\varphi$ is given in \eqref{eq:phi} and 
$u_0\approx 0.831443$ is the positive solution to
the equation $4\cosh^4u=\cosh 4u.$
\end{lem}

\begin{pf}
We observe that for $u>0,$
\begin{align*}
&\quad \varphi(2u)=2\,\mathrm{artanh}\,\big[(\sinh 2u)/\sqrt{\cosh 4u}\big]<2u \\
\Leftrightarrow &\quad
\frac{\sinh 2u}{\sqrt{\cosh 4u}}=\frac{2\sinh u\cosh u}{\sqrt{\cosh 4u}}
< \tanh u=\frac{\sinh u}{\cosh u} \\
\Leftrightarrow &\quad
4<\frac{\cosh 4u}{\cosh^4u}.
\end{align*}
Since $(\cosh 4u)/\cosh^4u$ increases from $1$ to $8$ when $u$ moves from
$0$ to $+\infty,$ there exists a unique number $u_0>0$ satisfying
the relation $4=(\cosh 4u_0)/\cosh^4u_0.$
We now see that $\varphi (2u)<2u$ if and only if $u>u_0.$
\end{pf}

The following elementary result is also needed later.

\begin{lem}\label{lem:sh}
The function $f(x)=x\slash\log(1+2\sinh x)$ strictly increases from $1/2$ to $1$
as $x$ moves from $0$ to $+\infty.$
\end{lem}

\begin{pf}
Because $f(x)=x/\log(e^x-e^{-x}+1),$ differentiation yields
\[f'(x)=h(x)/\big[\log(e^x-e^{-x}+1)\big]^2\,,\,\, \,{\rm  where}\,\,
h(x)=\log(e^x-e^{-x}+1)-\frac{x(e^x+e^{-x})}{e^x-e^{-x}+1}\,.
\]
Further, we have
$$
h'(x)=-\frac{x(e^x-e^{-x})}{e^x-e^{-x}+1}
+\frac{x(e^x+e^{-x})^2}{(e^x-e^{-x}+1)^2}
=\frac{x(e^{-x}-e^x+4)}{(e^x-e^{-x}+1)^2}
=\frac{2x(2-\sinh x)}{(1+2\sinh x)^2}\,.
$$
We now see that $h'(x)>0$ for $0<x<\,\mathrm{arsinh}\, 2$ and $h'(x)<0$ for $\,\mathrm{arsinh}\, 2<x.$
Since $h(0)=0$ and 
$$
h(x)=x+\log(1+e^{-x}-e^{-2x})-x\frac{1+e^{-2x}}{1+e^{-x}-e^{-2x}}
=O(xe^{-x})=o(1)
$$
as $x\to +\infty,$
the function $h(x)$ is positive for all $x>0.$
Hence, $f'(x)>0$ for all $x>0,$ which implies that $f(x)$ is strictly increasing in $x>0.$
It is easy to see that $f(x)\to 1/2$ as $x\to0$ and that $f(x)\to 1$ as $x\to+\infty.$
\end{pf}

We are ready to prove our result.

\begin{lem}\label{lem:extremal}
Let $u>0.$
Then the quantity $M(u)$ defined in \eqref{eq:M} is evaluated as
$$
M(u)=
\begin{cases}
\chi(u)
& ~\text{if}~ 0<u<u_0, \\
2\sinh 2u &~ \text{if}~ u_0\le u,
\end{cases}
$$
where $\chi(u)$ is given in \eqref{eq:chi} and $u_0
\approx 0.831443$ is the positive solution to
the equation $4\cosh^4u=\cosh 4u.$
Moreover, when $0<u<u_0,$ a set $E\in{\mathcal E}(u)$ satisfies
${d}(E)=M(u)$ if and only if $E=E^*(u).$
\end{lem}

\begin{pf}
We denote the circle $\partial\Delta(i,2u)$ by $C$ in the following.
Since every $E\in{\mathcal E}(u)$ is contained in the closed disk ${\overline{\Delta}}(i,2u),$
the diameter $d(E)$ is at most $2\sinh 2u$ by Lemma \ref{lem:circle}(ii).
Hence, we observe that 
\begin{equation*}
M(u)\le 2\sinh 2u,\quad u>0.
\end{equation*}

First, we assume that $u\ge u_0;$ equivalently by Lemma \ref{lem:u0},
$\varphi(2u)\le 2u.$
Let $z_1, z_2$ be the endpoints of the horizontal diameter of the boundary circle
$C=\partial\Delta(i,2u).$
Note that ${\,\operatorname{Im}\,} z_j=\cosh 2u>1.$
Then, by Lemma \ref{lem:circle}(iii), $h_{\mathbb H}(z_1,z_2)=\varphi(2u)\le 2u.$
Thus, $E=\{i,z_1,z_2\}\in{\mathcal E}(u)$ which implies ${d}(E)=2\sinh 2u\le M(u).$
Therefore, we have proved that $M(u)=2\sinh 2u.$
Note that the extremal set $E$ is not necessarily unique when $\varphi(2u)<2u$
(for instance, we can rotate the diameter a little about the Euclidean center of $C$).

Next, we assume that $u<u_0;$ namely, $\varphi(2u)>2u.$
We prove that there exists a set $E_0\in{\mathcal E}(u)$ 
attaining the supremum in \eqref{eq:M}; namely, $M(u)={d}(E_0).$
Indeed, by definition, we can find a sequence of sets $E_k$ in ${\mathcal E}(u)$
such that ${d}(E_k)\to M(u)$ as $k\to\infty.$
Since each $E\in{\mathcal E}(u)$ is contained in the closed hyperbolic disk ${\overline{\Delta}}(i,2u),$
by passing to a subsequence if necessary, we may assume
that $E_k=\{i, z_k, w_k\}$ and $z_k\to z_\infty$ and $w_k\to w_\infty$
as $k\to \infty$ for some $z_\infty, w_\infty\in{\overline{\Delta}}(i,2u).$
By continuity, we have ${d} (E_\infty)=M(u)$ for 
$E_\infty=\{i,z_\infty, w_\infty\}.$
We have to check that $E_\infty$ belongs to ${\mathcal E}(u).$
If $E_\infty$ consists only of two points, by Lemma \ref{lem:gm1},
$$\log(1+M(u))\le J_{\mathbb H}(E_\infty)\le h_{\mathbb H}(E_\infty)=2u,$$
which contradicts Lemma \ref{lem:11/4} because $u\le u_0<\log(11/4).$
We have proved the claim.

Now assume that $E_0=\{i,z_0,w_0\}\in{\mathcal E}(u)$ satisfies ${d}(E_0)=M(u).$
By assumption, we have $z_0\in{\overline{\Delta}}(i, 2u)\cap {\overline{\Delta}}(w_0, 2u).$
Observe that $z_0\in\partial\Delta(i,2u)=C$ in the present situation.
In fact, let $r=h_{\mathbb H}(z_0,w_0)$ and suppose $h_{\mathbb H}(z_0,i)<2u.$
Then $z_0$ can be moved along the circle $\partial\Delta(w_0,r)$ upwards a bit
to get a new point $z_0'$ in such a way that 
$${\,\operatorname{Im}\,} z_0<{\,\operatorname{Im}\,} z_0',~ h_{\mathbb H}(z_0',i)<2u,~
h_{\mathbb H}(z_0',w_0)=r \quad {\rm and}\,\, |z_0-w_0|<|z_0'-w_0|$$ 
by Lemma \ref{lem:circle} and Remark \ref{rem:circle}.
Hence we would have $h_{\mathbb H}(E_0')=h_{\mathbb H}(E_0)$ and
${d}(E_0)<{d}(E_0')$ for $E_0'=\{i,z_0',w_0\}.$
This, however, violates the initial assumption that ${d}(E_0)=M(u).$
Therefore, we conclude that $h_{\mathbb H}(z_0,i)=2u.$
In the same way, we obtain $h_{\mathbb H}(w_0,i)=2u.$
We can further prove, as before (cf.~the proof of Lemma \ref{lem:11/4}), 
that $|z_0-w_0|={d}(E_0).$

The remaining task is now to determine the configuration
of the points $z_0, w_0$ on the circle $C$
maximizing the quantity $|z_0-w_0|$ under the constraints
$h_{\mathbb H}(z_0,w_0)\le 2u$ and $\min\{{\,\operatorname{Im}\,} z_0, {\,\operatorname{Im}\,} w_0\}\ge 1.$
We recall that the hyperbolic distance of the endpoints of an arbitrary 
Euclidean diameter of $C$ is at least $\varphi(2u)$ by Lemma \ref{lem:circle}(iii).
We first suppose that $\varphi(2u)<2u.$
Let $C_0$ be the shorter component of $C\setminus\{z_0, w_0\}.$
It is evident that the chord $|z_0-w_0|$ is shortest when (and only when)
$z_0$ and $w_0$ are situated symmetrically with respect to the imaginary axis.
Therefore, we have 
\[ E_0=E^*(u)\quad {\rm and}\,\, M(u)=2u/\log(1+{d}(E^*(u)))=\xi(u)\,.\]

By the above proof, uniqueness of the extremal set for $0<u\le u_0$ is clear.
Thus, the proof is now complete.
\end{pf}

\begin{rem}\label{rem:M}
In view of Lemmas \ref{lem:11/4} and \ref{lem:sh}, as a corollary of the last lemma,
we have the inequality
$$
\inf_{E\in{\mathcal E}(u)} \frac{2u}{J_{\mathbb H}(E)}
=\frac{2u}{\log(1+M(u))}<1
$$
for every $u>0.$
\end{rem}

We are now in a position to prove the following theorem.

\begin{thm}\label{thm:het}
There is a zero $u=u^*$ of the derivative $\xi'(u)$ of the function
$$
\xi(u)=\frac{2u}{\log(1+\chi(u))}
$$
in the interval $0<u<u_0\approx 0.83$ such that
$$
\kappa({\mathbb H})=\frac{h_{\mathbb H}(z^*,w^*)}{\log(1+|z^*-w^*|)}
=\frac{h_{\mathbb H}(E^*)}{\log(1+{d}(E^*)/d(E^*,\partial{\mathbb H}))},
$$
where $u_0$ is given in Lemma \ref{lem:u0},
$E^*=E^*(u^*)=\{i,z^*, w^*\},$
$z^*=ie^{t^*+i\theta^*},~ w^*=ie^{t^*-i\theta^*}$
and $t^*, \theta^*$ are given in \eqref{eq:t} and \eqref{eq:theta}, respectively,
for $u=u^*.$
Moreover, if $\kappa({\mathbb H})=h_{\mathbb H}(E)/\log(1+{d}(E)/d(E,\partial{\mathbb H}))$
for a three-point set $E$ in ${\mathbb H},$ then there are real numbers $a,b$
with $a>0$ such that $E=aE^*+b.$
\end{thm}

\begin{pf}
Lemma \ref{lem:extremal} implies that for $u\ge u_0=0.831\dots,$
$$
\frac{2u}{\log(1+M(u))}=\frac{2u}{\log(1+2\sinh 2u)}.
$$
Since the function $x/\log(1+2\sinh x)$ is increasing in $0<x<+\infty$
by Lemma \ref{lem:sh}, we can restrict the range of the infimum in \eqref{eq:g3}
to $(0,u_0]$:
$$
\kappa({\mathbb H})=\kappa_3({\mathbb H})
=\inf_{0<u\le u_0}\frac{2u}{\log(1+M(u))}
=\inf_{0<u\le u_0}\frac{2u}{\log(1+\chi(u))}
=\inf_{0<u\le u_0}\xi(u),
$$
where $\chi(u)$ is given in \eqref{eq:chi}.
By the form of $\chi(u)$ in \eqref{eq:chi}, we observe that
$\chi(u)=2u+2\sqrt3 u^2+O(u^3)$ as $u\to0^+.$
Thus, we obtain $\xi(u)\ge 2u/\log(1+\chi(u))=1-(\sqrt3-1)u+O(u^2)$ as $u\to0^+.$
In particular, $\xi(0^+)=1$ and $\xi'(0^+)=1-\sqrt 3<0.$
Since $\xi'(u_0)=0.1917\dots>0,$
the above infimum of $\xi(u)$ is attained at its critical point in $(0,u_0).$

The last assertion easily follows from the uniqueness of the extremal set in
Lemma \ref{lem:extremal}.
The proof is now complete.
\end{pf}

See Figure 2 for the graph of the function $2u/\log(1+M(u)).$
By numerical computations, we obtain
$u^*\approx       0.432335123777,
t^*\approx          0.727535978839,
\theta^*\approx 0.419463976058,$
and $\kappa({\mathbb H})=
\xi(u^*)\approx 0.8750987500145.$
Note that by Theorem \ref{thm:tmt} $\kappa({\mathbb H})=\kappa(\Omega)$ for a convex hyperbolic domain $\Omega$. In conclusion, we have the following theorem:

\begin{thm}\label{thm_approxforOmega}
For any convex hyperbolic domain $\Omega$, $\kappa(\Omega)\approx0.875098750014$.
\end{thm}

\begin{figure}[htbp]
\begin{center}
\includegraphics[width=0.6\textwidth]{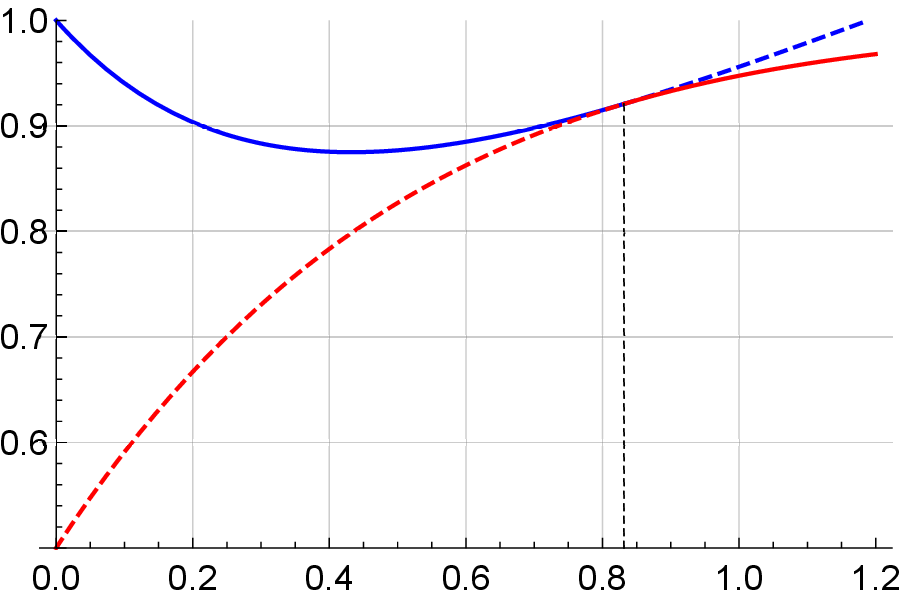}
\caption{The graph of $2u/\log(1+M(u))$ (the thick line); 
the blue curve indicates the graph of $\xi(u)$ and 
the red one does the graph of $2u/\log(1+2\sinh 2u)$}
\end{center}
\end{figure}

Finally, we prove Theorem \ref{thm:existence}.

\begin{pf}[Proof of Theorem \ref{thm:existence}]
It remains to prove the first assertion.
Let $\Omega\subsetneq{\mathbb C}$ be a convex domain and suppose that
$\kappa(\Omega)=h_\Omega(E)/J_\Omega(E)$ for a compact
subset $E$ of $\Omega.$
As in the proof of Lemma \ref{lem:gamma_n}, we 
take points $z_0, z_1, z_2\in E$ so that ${d}(E)=|z_1-z_2|$ 
and $d(E, \partial \Omega)=d(z_0, \partial \Omega)$
and let $E_0=\{z_0, z_1, z_2\}.$
(Since $\kappa(\Omega)=\kappa({\mathbb H})<1,$ the set $E_0$ contains exactly three
points.)
By Lemma \ref{lem:gamma_n}, we have $\kappa_3(\Omega)=\kappa(\Omega).$
Thus, in the chain of inequalities \eqref{eq:chain}, the last term is the same as
the initial term.
Thus, we have $h_\Omega(E)=h_\Omega(E_0).$
Hence $\kappa(\Omega)=h_\Omega(E_0)/J_\Omega(E_0).$

Let $\zeta_0\in\partial\Omega$ be such that $d(E_0,\partial\Omega)
=d(z_0,\partial\Omega)=|z_0-\zeta_0|.$
Take a half-plane $H$ as in the proof of Lemma \ref{lem:conv}
such that $\Omega\subset H$ and $z_0\in\partial H.$
Then $J_\Omega(E_0)=J_H(E_0)$ and $h_\Omega(E_0)\ge h_H(E_0).$
If $\Omega$ is a proper subdomain of $H,$ then we would have
$h_H(E_0)<h_\Omega(E_0).$
Thus,
$$
\kappa(H)\le\frac{h_H(E_0)}{J_H(E_0)}
<\frac{h_\Omega(E_0)}{J_\Omega(E_0)}=\kappa(\Omega).
$$
On the other hand, Theorem \ref{thm:tmt} yields $\kappa(H)=\kappa(\Omega),$
which is a contradiction.
Thus, $\Omega$ equals $H,$ a half-plane.
\end{pf}


\section{Application to Capacity Estimation}


Finally, we apply the results above to capacity estimation. 
First, we recall some basic notions.

\begin{defn}\cite[Def. 9.2, p. 150]{HKV20}
A pair $(\Omega,E)$ of a domain $\Omega$ in ${\mathbb C}$ and a non-empty compact subset 
$E$ of $\Omega$ is called a \emph{condenser}. 
The \emph{capacity} of this condenser is defined to be
\begin{align*}
\mathrm{cap}\,(\Omega,E)=\inf_u\iint_{{\mathbb C}}|\nabla u(z)|^2 dxdy \quad (z=x+iy),
\end{align*}
where the infimum is taken over the family of all non-negative functions $u$ in the Sobolev
class $W_{\rm{loc}}^{1,2}({\mathbb C})$ with compact support in $\Omega$ such that $u(z)\geq1$ for $z\in E$.
\end{defn}

If $\Omega$ is a simply connected proper subdomain of ${\mathbb C}$
and $E$ is a (non-degenerate) continuum in $\Omega$ such that the set $R=\Omega\setminus E$ is
a doubly connected domain (a ring), then its modulus is known to be $2\pi/\mathrm{cap}\,(\Omega,E).$

We define the homeomorphism $\mu:(0,1)\to{\mathbb R}^+$ by the formula (see, e.g., \cite[7.4.1, p. 122]{HKV20})
$$
\mu(r)=\frac{\pi}{2}\cdot\frac{\mathcal{K}(\sqrt{1-r^2})}{\mathcal{K}(r)},
$$
where $\mathcal{K}(r)$ is Legendre's complete elliptic integral of the first kind defined by
$$
\mathcal{K}(r)=\int_0^1\frac{dx}{\sqrt{(1-x^2)(1-r^2x^2)}}.
$$
It is known that $\mu(r)$ represents the modulus of the Gr\"otzsch ring ${\mathbb D}\setminus [0,r].$
In particular, $\mu(r)$ decreases from $+\infty$ to $0$ as $r$ moves from $0$ to $1.$
We note that $2\pi/\mu(r)$ is the capacity of ${\mathbb D}\setminus [0,r].$
For later convenience, we put
$$
\Phi(x)=\frac{2\pi}{\mu(\tanh(x/2))},\quad 0<x<\infty.
$$
Note that $\Phi(x)$ increases from $0$ to $+\infty$ as $x$ moves from $0$ to $+\infty.$
We are ready to give the main result in this section.
Recall that $J_\Omega(E)=\log(1+d(E)/{d}(E,\partial\Omega)).$

\begin{thm}\label{thm:cap-lower}
Let $E$ be a continuum in a simply connected domain $\Omega\subsetneq{\mathbb C}.$
Then the following are valid.
\begin{enumerate}
\item[(i)]
The inequality
$$
\mathrm{cap}\,(\Omega, E)
\ge \Phi(\kappa(\Omega)J_\Omega(E))
\ge \Phi(\kappa_0J_\Omega(E))
$$
holds, where $\kappa_0$ is given in \eqref{eq:sc}.
\item[(ii)]
If $\Omega$ is convex,
$$
\mathrm{cap}\,(\Omega, E) \ge  \Phi(\kappa_1 J_\Omega(E))),
$$
where $\kappa_1=\kappa({\mathbb D})>0.87509875.$
\end{enumerate}
\end{thm}

\begin{pf}
Let $f:\Omega\to{\mathbb D}$ be a conformal homeomorphism and set $E'=f(E).$
Since the capacity and the hyperbolic distance are conformally invariant, we obtain
$$
\mathrm{cap}\,(\Omega,E)=\mathrm{cap}\,({\mathbb D},E')\ge \Phi(h_{\mathbb D}(E'))=\Phi(h_\Omega(E)),
$$
where we used a consequence of the circular symmetrization (see \cite[Lemma 9.20, p. 163]{HKV20}).
Other parts follow from Corollary \ref{cor} and Theorem \ref{thm_approxforOmega}.
\end{pf}

\begin{ex}
Consider next an example where $\Omega=\{z: -1<{\,\operatorname{Im}\,} z<1\}$ and $E=[1,2]$. 
Because $\Omega$ is convex, it follows from Theorem \ref{thm:cap-lower} that
\begin{align*}
\mathrm{cap}\,(\Omega,E)
\geq \Phi(\kappa_1 J_\Omega(E))
\approx\frac{2\pi}{\mu(0.43754937\log 2))}>2.4288.
\end{align*}
By applying the circular (spherical) symmetrization 
(see \cite[9.1, pp. 155-157]{HKV20}) 
with the origin as a center and $x$-axis as the symmetrization axis. 
Observe first that
the negative $x$-axis is contained in the complement of the symmetrized
condenser whereas $[1,2]$ remains invariant and hence
\begin{align*}
\mathrm{cap}\,(\Omega,E)\geq\tau_2(1)=2,  
\end{align*}
where $\tau_2(t)$ denotes the capacity of the Teichm\"uller ring
${\mathbb C}\setminus([-1,0]\cup [t,+\infty))$ for $t>0$
(see \cite[7.3, pp. 120]{HKV20}),
which is  a weaker lower bound for the capacity than what we proved above.
On the other hand, if we take into account that the whole left half-plane
is contained in the complement of the symmetrized condenser, we obtain
\begin{align*}
\mathrm{cap}\,(\Omega,E)\geq \Phi(\log 2)=
\frac{2 \pi}{\mu(\tanh(\log\sqrt 2))}\approx 2.55852.
\end{align*}
Hence the value of our bound given in Theorem \ref{thm:cap-lower}
lies between these two bounds obtained by symmetrization.
Finally, let us find the exact value of $\mathrm{cap}\,(\Omega,E).$
Obviously, $\mathrm{cap}\,(\Omega,E)=\mathrm{cap}\,(\Omega,E_0),$ where $E_0=[0,1].$
Note that the function $\displaystyle f(z)=\frac2\pi\log\frac{1+z}{1-z}$ maps
the unit disk ${\mathbb D}$ onto $\Omega$ and that $f^{-1}(E_0)=[0,\tanh(\pi/4)].$
Thus
$$
\mathrm{cap}\,(\Omega,E)=\mathrm{cap}\,(\Omega,E_0)=\mathrm{cap}\,({\mathbb D},[0,\tanh(\pi/4)])
=\frac{2\pi}{\mu(\tanh(\pi/4))}=\Phi(\frac{\pi}{2})\approx 3.75108.
$$
\end{ex}


{\bf Acknowledgements.} The authors would like to thank the referee for detailed and constructive corrections.


\def\cprime{$'$} \def\cprime{$'$} \def\cprime{$'$}
\providecommand{\bysame}{\leavevmode\hbox to3em{\hrulefill}\thinspace}
\providecommand{\MR}{\relax\ifhmode\unskip\space\fi MR }
\providecommand{\MRhref}[2]{%
  \href{http://www.ams.org/mathscinet-getitem?mr=#1}{#2}
}
\providecommand{\href}[2]{#2}

\end{document}